\documentclass[12pt]{article}
\newtheorem{lemma}{Lemma}[section]
\newtheorem{theorem}[lemma]{Theorem}
\newtheorem{proposition}[lemma]{Proposition}

\newtheorem{definition}[lemma]{Definition}
\newenvironment{proof}{{\bf Proof}}{{\hfill $ \Box $}\vskip 4mm}
\newenvironment{remark}{\addtocounter{lemma}{1}
{\bf Remark \thelemma}}{{\hfill}\vskip 4mm}
\newenvironment{remarks}{\addtocounter{lemma}{1}
{\bf Remarks \thelemma}}{{\hfill}\vskip 4mm}

\newenvironment{examples}{\addtocounter{lemma}{1}
{\bf Examples \thelemma}}{{\hfill}\vskip 4mm}
\newcommand{\nc}{\newcommand}
\nc{\rnc}{\renewcommand}
\nc{\nt}{\newtheorem}

        {~\hfill~\fbox{}\end{trivlist}}

\nc{\thlabel}[1]{\label{theo:#1}}
\nc{\thref}[1]{Theorem~\ref{theo:#1}}
\nc{\selabel}[1]{\label{sect:#1}}
\nc{\seref}[1]{Section~\ref{sect:#1}}
\nc{\lelabel}[1]{\label{lemm:#1}}
\nc{\leref}[1]{Lemma~\ref{lemm:#1}}
\nc{\prlabel}[1]{\label{prop:#1}}
\nc{\prref}[1]{Proposition~\ref{prop:#1}}
\nc{\colabel}[1]{\label{coro:#1}}
\nc{\coref}[1]{Corollary~\ref{coro:#1}}
\nc{\exlabel}[1]{\label{exam:#1}}
\nc{\exref}[1]{Example~\ref{exam:#1}}
\nc{\delabel}[1]{\label{defi:#1}}
\nc{\deref}[1]{Definition~\ref{defi:#1}}
\nc{\eqlabel}[1]{\label{equation:#1}}
\nc{\eqref}[1]{(\ref{equation:#1})}
\nc{\csm}{\mbox{$\triangleright\!\!\!<$}}
\nc{\smc}{\mbox{$>\!\!\!\triangleleft$}}
\nc{\trr}{\triangleright}
\providecommand{\operatorname}[1]{\mathrm{#1}\,}
\nc{\Hom}{\operatorname{Hom}}
\nc{\Mor}{\operatorname{Mor}}
\nc{\Aut}{\operatorname{Aut}}
\nc{\Ann}{\operatorname{Ann}}
\nc{\Ker}{\operatorname{Ker}}
\nc{\Trace}{\operatorname{Trace}}
\nc{\Char}{\operatorname{Char}}
\nc{\Mod}{\operatorname{Mod}}
\nc{\End}{\operatorname{End}}
\nc{\Spec}{\operatorname{Spec}}
\nc{\Span}{\operatorname{Span}}
\nc{\sgn}{\operatorname{sgn}}
\nc{\Id}{\operatorname{Id}}
\nc{\Com}{\operatorname{Com}}

\def\Box{\mbox{$\sqcap\!\!\!\!\sqcup$}}

\nc{\dht}{\mbox{$\rightharpoonup\hspace{-2ex}\rightharpoonup$}}
\nc{\dhtb}{\mbox{$\leftharpoonup\hspace{-2ex}\leftharpoonup$}}
\nc{\nd}{\mbox{$\not|$}} 
\providecommand{\text}[1]{\mbox{{\textrm #1}}}

\nc{\nci}{\mbox{$\not\subseteq$}}
\nc{\scontainin}{\mbox{$\mbox{}\subseteq\hspace{-1.5ex}\raisebox{-.5ex}{$_\prime
$}\hspace*{1.5ex}$}}

\def\ot{\otimes}

\def\doublerightleft#1#2{{\lower.2ex\vbox{
\hbox{${\smash{\mathop{\longrightarrow}\limits^{#1}}}$}\vspace*{-4mm}
\hbox{${\smash{\mathop{\longleftarrow}\limits_{#2}}}$}}}}

\def\nint{\hbox{\bbb Z}}
\def\nrat{\hbox{\bbb Q}}

\newfont{\bbb}{msbm10 scaled\magstep1}  
\newfont{\bbbsub}{msbm10}                
\newfont{\msam}{msam10 scaled\magstep1}

\begin{document}
\title{The Long dimodules category and nonlinear equations}
\author{G. Militaru
\\University of Bucharest
\\Faculty of Mathematics\\Str. Academiei 14
\\RO-70109 Bucharest 1, Romania
\\e-mail: gmilit@al.math.unibuc.ro}
\date{}
\maketitle
\begin{abstract}
\noindent Let $H$ be a bialgebra and ${}_H{\cal L}^H$ be the category of
Long $H$-dimodules defined, for a commutative and cocommutative $H$, by F. W.
Long in \cite{Lo} and studied in connection to the Brauer group of a so
called $H$-dimodule algebra. For a commutative and cocommutative $H$,
${}_H{\cal L}^H={}_H{\cal YD}^H$ (the category of
Yetter-Drinfel'd modules), but for an arbitrary $H$ the categories
${}_H{\cal L}^H$ and ${}_H{\cal YD}^H$ are basically different.
Keeping in mind that the category ${}_H{\cal YD}^H$ is deeply involved in
solving the quantum Yang-Baxter equation, we shall study the category
${}_H{\cal L}^H$ of $H$-dimodules in connection with what we have called
the ${\cal D}$-equation: $R^{12}R^{23}=R^{23}R^{12}$, where
$R\in \End_k(M\ot M)$ for a vector space $M$ over a field $k$. The main
result is a FRT type theorem: if $M$ is finite dimensional, then
any solution $R$ of the ${\cal D}$-equation has the form
$R=R_{(M,\cdot ,\rho )}$, where $(M,\cdot,\rho)$ is a Long $D(R)$-dimodule
over a bialgebra $D(R)$ and $R_{(M,\cdot ,\rho )}$ is the special map
$R_{(M,\cdot,\rho)}(m\ot n):=\sum n_{<1>}\cdot m\ot n_{<0>}$.
In the last section, if $C$ is a coalgebra and $I$ is a coideal of $C$, we shall
introduce the notion of ${\cal D}$-map on $C$, that is a $k$-bilinear
map $\sigma :C\ot C/I\to k$ satisfying a condition which ensures on one hand
that for any right $C$-comodule, the special map $R_{\sigma}$ is a solution
of the ${\cal D}$-equation, and on the other that, in the finite case,
any solution of the ${\cal D}$-equation  has this form.
\end{abstract}

\section{Introduction}
Let $k$ be a field, $M$ a $k$-vector space, $R\in \End_k(M\ot M)$
such that a certain equation $({\cal E})$ holds in $\End_k(M\ot M\ot M)$.
The starting point of the present paper is the following general question:

{\em "Does the Hopf algebra theory offer an effective technique for solving
the equation $({\cal E})$"?}

Using as a source of inspiration the quantum Yang-Baxter equation, in the
solution of which the instruments offered by Hopf algebra theory have proven
to be very efficient (see \cite{FRT}, \cite{K}, \cite{LR}, \cite{R1}),
we can identify two ways of approaching the above problem:

1. Let $H$ be a bialgebra. The first idea is to define a new category
${}_H{\cal E}^H$: the objects in this category are threetuples
$(M,\cdot ,\rho)$, where $(M,\cdot)$ is a left $H$-module and
$(M,\rho)$ is a right $H$-comodule such that $(M,\cdot ,\rho)$
satisfies a compatibility relation which ensures that the special map
$$
R_{(M,\cdot,\rho)}:M\ot M\to M\ot M,\quad
R_{(M,\cdot,\rho)}(m\ot n):=\sum n_{<1>}\cdot m\ot n_{<0>}
$$
for all $m$, $n\in M$, is a solution of the equation $({\cal E})$. The deal
is now the converse: if $M$ is a finite dimensional vector space and $R$
is a solution of the equation $({\cal E})$, can we construct a bialgebra
$E(R)$ such that $(M,\cdot ,\rho)$ has a structure of object in
${}_{E(R)}{\cal E}^{E(R)}$ and $R=R_{(M,\cdot ,\rho )}$? Roughly speaking,
an affirmative answer to this FRT type theorem tells us that the maps
$R_{(M,\cdot ,\rho )}$, where $(M,\cdot ,\rho)$ is an object in
${}_H{\cal E}^H$ for some bialgebra $H$, "count all solutions" of the
equation $({\cal E})$. Two questions arise in a natural way: first, which is
the bialgebra $E(R)$?
Suppose that $M$ has dimension $n$. Let ${\cal M}^n(k)$ be the comatrix
coalgebra of order $n$ and $(T({\cal M}^n(k)),m,1,\Delta, \varepsilon)$ the
unique bialgebra structure of the tensor algebra $(T({\cal M}^n(k)),m,1)$
whose comultiplication $\Delta$ and counit $\varepsilon$ extend the
comultiplication and the counit from ${\cal M}^n(k)$.
As a general rule, $E(R)=T({\cal M}^n(k))/I$, where $I$ is a bi-ideal of
$T({\cal M}^n(k))$ generated by some "obstructions".
The second question arising here regards the potential benefits of such an
approach to the nonlinear equation $({\cal E})$. The answer to this question
is easy: following the technique evidenced by Radford in \cite{R2} for the
quantum Yang-Baxter equation, reduced to the pointed case, we can classify
the solutions of the equation $({\cal E})$, at least for low dimensions.

2. The second approach to our problem is the following: let $H$ be a
bialgebra and $\sigma:H\ot H\to k$ a $k$-bilinear map.
We require for the map $\sigma$ to satisfy certain properties $(P)$ which
ensure that for any right $H$-comodule $M$, the natural map
$$
R_{\sigma}:M\ot M\to M\ot M,\quad
R_{\sigma}(m\ot n)=\sum \sigma(m_{<1>}\ot n_{<1>})m_{<0>}\ot n_{<0>}
$$
for all $m$, $n\in M$, is a solution for the equation $({\cal E})$.
Similar to case 1, the converse is interesting: in the finite dimensional
case, any solution arises in this way.

In the case of the quantum Yang-Baxter equation the role of the category
${}_H{\cal E}^H$ is played by ${}_H{\cal YD}^H$, the category of
Yetter-Drinfel'd modules (or crossed modules) defined by Yetter in
\cite{Y}. Of course, here we have in mind the version of the famous FRT
theorem given by Radford in \cite{R1}. The role of the maps $\sigma$ which
satisfy the properties $(P)$ is played by the co-quasitriangular (or braided)
bialgebras. For more information regarding the Yetter-Drinfel'd modules or
co-quasitriangular bialgebras we refer to \cite{K}, \cite{R1},
\cite{R2}, or to the more recent \cite{CMZ1}, \cite{CMZ2}.

Recentely, in \cite{M1} we evidence the fact that the clasical
category ${}_H{\cal M}^H$ of $H$-Hopf modules (\cite{A}) is also
deeply involved in solving a certain non-linear equation.
We called it the Hopf equation, and it is:
$$
R^{12}R^{23}=R^{23}R^{13}R^{12}
$$
The Hopf equation is equivalent, taking $W=\tau R\tau$, where $\tau$ is
the flip map, with the pentagonal equation
$$
W^{12}W^{13}W^{23}=W^{23}W^{12}
$$
which plays a fundamental role in the duality theory of operator algebras
(see \cite{BS} and the references indicated here). We remember that the
unitary fundamental operator defined by Takesaki for a Hopf von Neumann
algebra is a solution of the Hopf equation (see lemma 4.9 of \cite{Ta}).

Although at first sight the categories ${}_H{\cal M}^H$ and
${}_H{\cal YD}^H$ are completely different, in fact they are particular
cases of the same general category ${}_A{\cal M}(H)^C$ of
Doi-Hopf modules defined by Doi in \cite{D}: if $A$ is a right
$H$-comodule algebra and $C$ is a left $H$-module coalgebra then a
threetuple $(M,\cdot, \rho)$ is an object in ${}_A{\cal M}(H)^C$ if
$(M,\cdot)$ is a left $A$-module, $(M,\rho)$ is a right $C$-comodule and
the following compatibility condition holds
\begin{equation}\label{cin}
\rho(a\cdot m)=\sum a_{<0>}\cdot m_{<0>}\ot a_{<1>}m_{<1>}
\end{equation}
for all $a\in A$, $m\in M$. Taking $A=C=H$, then
${}_H{\cal M}(H)^H={}_H{\cal M}^H$, the category of classical $H$-Hopf
modules; on the other hand ${}_H{\cal YD}^H={}_H{\cal M}(H^{op}\ot H)^H$,
where $H$ can be viewed (see \cite{CMZ1}) as an $H^{op}\ot H$-module
(comodule) coalgebra (algebra).

Now, let $A=H$ with the right $H$-comodule structure via $\Delta$ and
$C=H$ with the trivial $H$-module structure:
$h\cdot k:=\varepsilon(h)k$ for all $h$, $k\in H$. Then the compatibility
condition (\ref{cin}) takes the form
$$
\rho(h\cdot m)=\sum h\cdot m_{<0>}\ot m_{<1>}
$$
for all $h\in H$, $m\in M$. We shall denote this category with
${}_H{\cal L}^H$, which is also a special case of the Doi-Hopf module
category. This category ${}_H{\cal L}^H$ was defined first by F.W. Long in
\cite{Lo} for a commutative and cocommutative $H$ and was studied in
connection to the construction of the Brauer group of an $H$-dimodule
algebra (an object in ${}_H{\cal L}^H$ was called in \cite{Lo} an
$H$-dimodule). It is interesting to note that for a commutative and
cocommutative $H$ the category ${}_H{\cal YD}^H$ of Yetter-Drinfel'd
modules is precisely ${}_H{\cal L}^H$, the category of Long $H$-dimodules
(see, for example, \cite{CMZ1}). Of course, for an arbitrary $H$, the
categories ${}_H{\cal YD}^H$ and ${}_H{\cal L}^H$ are basically different.
Keeping in mind that the category ${}_H{\cal YD}^H$ of Yetter-Drinfel'd
modules plays a determinant role in describing the solutions of the quantum
Yang-Baxter equation, the following question is natural:

{\sl "In which equation will the category ${}_H{\cal L}^H$ of
Long $H$-dimodules play a key role ?"}

Answering this question would complete the image of the way in which
the three well known categories from Hopf algebras theory,
${}_H{\cal M}^H$, ${}_H{\cal YD}^H$ and ${}_H{\cal L}^H$,
participate in solving nonlinear equations.

In this paper we shall study in detail what we have called the
${\cal D}$-equation: $R\in \End_k(M\ot M)$ is called a solution of
the ${\cal D}$-equation if
$$
R^{12}R^{23}=R^{23}R^{12}.
$$
In the study of this equation, we shall apply the general treatment
presented above.
In theorem 3.5 we prove that in the finite dimensional case any solution $R$
of the ${\cal D}$-equation has the form $R=R_{(M,\cdot ,\rho )}$, where
$(M,\cdot ,\rho)\in {}_{D(R)}{\cal L}^{D(R)}$, for some bialgebra $D(R)$.
The bialgebra $D(R)$ is a quotient of $T({\cal M}^n(k))$ by a bi-ideal
generated by some obstruction elements $o(i,j,k,l)$ which have degree one
in the graded algebra $T({\cal M}^n(k))$, i.e. $o(i,j,k,l)$ lie in the
coalgebra ${\cal M}^n(k)$.
This observation led us in the last section to define the concept of
${\cal D}$-map in a more general case, relative only to a coalgebra.
The ${\cal D}$-equation is obtained from the quantum Yang-Baxter equation
by deleting the middle term from both sides.
This operation destroys a certain symmetry which exists in the case of the
quantum Yang-Baxter equation, and it is reflected in the fact that the
${\cal D}$-maps are defined as $k$-bilinear maps from $C\ot C/I$ to $k$,
where $I$ is a coideal of a coalgebra $C$. This "assymetry" of the
${\cal D}$-equation is also underlined in the second item of the main
theorem of this section. As applications, we present several examples
of such constructions.

Following Doi's philosophy (\cite{D}), we can define a more general
category ${}_{A}{\cal L}^{C}$, where $A$ is an algebra and $C$ is a colagebra.
Any decomposable left $A$-module can be viewed as an object of
${}_{A}{\cal L}^{C}$, where $C=k[X]$, for some set $X$.
We can thus view the category  ${}_{A}{\cal L}^{C}$ as a generalisation of
the decomposable modules category.

{\sl Acknowledgement}: the author thanks S. Caenepeel for his helpful
comments on a preliminary version of this paper.

\section{Preliminaries}
Throughout this paper, $k$ will be a field.
All vector spaces, algebras, coalgebras and bialgebras that we consider
are over $k$. $\ot$ and $\Hom$ will mean $\ot_k$ and $\Hom_k$.
For a coalgebra $C$, we will use Sweedler's $\Sigma$-notation, that is,
$\Delta(c)=\sum c_{(1)}\ot c_{(2)},~(I\ot\Delta)\Delta(c)=
\sum c_{(1)}\ot c_{(2)}\ot c_{(3)}$, etc. We will
also use  Sweedler's notation for right $C$-comodules:
$\rho_M(m)=\sum m_{<0>}\otimes m_{<1>}$, for any $m\in M$ if
$(M,\rho_M)$ is a right $ C$-comodule. ${\cal M}^C$ will be the
category of right $C$-comodules and $C$-colinear maps and
${}_A{\cal M}$ will be the category of left $A$-modules and
$A$-linear maps, if $A$ is a $k$-algebra. An important role in the present
paper will be played by the ${\cal M}^n(k)$, the comatrix coalgebra of
order $n$, i.e. ${\cal M}^n(k)$ is the $n^2$-dimensional vector space
with $\{c_{ij}\mid i,j=1,\cdots,n \}$ a $k$-basis such that
\begin{equation}\label{com1}
\Delta(c_{jk})=\sum_{u=1}^{n}c_{ju}\ot c_{uk},
\quad \varepsilon (c_{jk})=\delta_{jk}
\end{equation}
for all $j,k=1,\cdots, n$. We view $T({\cal M}^n(k))$ with the unique
bialgebra structure which can be defined on the tensor algebra
$T({\cal M}^n(k))$ which extend the comultiplication $\Delta$ and the
counity $\varepsilon$ of ${\cal M}^n(k)$.

For a vector space $M$, $\tau :M\otimes M\to M\otimes M$
will denote the flip map, that is, $\tau (m\otimes n)=n\otimes m$
for all $m$, $n \in M$. $\tau^{(123)}$ will be the automorphism of
$M\ot M\ot M$ given by $\tau^{(123)}(l\ot m\ot n)=n\ot l\ot m$, for all
$l$, $m$, $n \in M$. If $R:M\ot M\to M\ot M$ is a linear map
we denote by $R^{12}$, $R^{13}$, $R^{23}$ the maps of $\End_k(M\ot M\ot M)$
given by
$$
R^{12}=R\ot I, \quad R^{23}=I\ot R,\quad
R^{13}=(I\ot \tau)(R\ot I)(I\ot \tau).
$$
\section{The ${\cal D}$-equation}
We will introduce the following

\begin{definition}
Let $M$ be a vector space and $R\in \End_k(M\ot M)$. We shall say
that $R$ is a solution for the ${\cal D}$-equation if
\begin{equation}\label{Ceq}
R^{12}R^{23}=R^{23}R^{12}
\end{equation}
in $\End_k(M\ot M\ot M)$.
\end{definition}

Solving the ${\cal D}$-equation is not an easy task. Even for the two-dimensional
case, solving the ${\cal D}$-equation is equivalent to solving a homogenous
system of 64 nonlinear equations (see equation (\ref{acum})).

\begin{proposition}\label{hei}
Let $M$ be a finite dimensional vector space and $\{m_1,\cdots,m_n \}$
a basis of $M$. Let $R$, $S\in \End_k(M\ot M)$ given by
$$
R(m_v\ot m_u)=\sum_{i,j}x_{uv}^{ji}m_i\ot m_j, \quad
S(m_v\ot m_u)=\sum_{i,j}y_{uv}^{ji}m_i\ot m_j,
$$
for all $u$, $v=1,\cdots ,n$, where $(x_{uv}^{ji})_{i,j,u,v}$,
$(y_{uv}^{ji})_{i,j,u,v}$ are two families of scalars of $k$. Then
$$R^{23}S^{12}=S^{12}R^{23}$$
if and only if
\begin{equation}\label{unde}
\sum_{v}x_{kv}^{ji}y_{lq}^{vp}=
\sum_{\alpha}x_{kl}^{j\alpha}y_{\alpha q}^{ip}
\end{equation}
for all $i$, $j$, $k$, $l$, $p$, $q=1,\cdots, n$. In particular, $R$ is a
solution for the C-equation if and only if
\begin{equation}\label{acum}
\sum_{v}x_{kv}^{ji}x_{lq}^{vp}=
\sum_{\alpha}x_{kl}^{j\alpha}x_{\alpha q}^{ip}
\end{equation}
for all $i$, $j$, $k$, $l$, $p$, $q=1,\cdots, n$.
\end{proposition}

\begin{proof}
For $k$, $l$, $q=1,\cdots, n$ we have:
\begin{eqnarray*}
R^{23}S^{12}(m_q\ot m_l\ot m_k )&=&
R^{23}\Bigl(\sum_{v,p}
y_{lq}^{vp}m_{p}\ot m_v\ot m_{k} \Bigl)\\
&=&\sum_{v,p,i,j}x_{kv}^{ji}y_{lq}^{vp}m_{p}\ot m_i\ot m_j \\
&=& \sum_{i,j,p} \Bigl ( \sum_{v}x_{kv}^{ji}
y_{lq}^{vp}\Bigl ) m_{p}\ot m_i\ot m_j
\end{eqnarray*}
and
\begin{eqnarray*}
S^{12}R^{23}(m_q\ot m_l\ot m_k )&=&
S^{12}\Bigl(\sum_{j,\alpha}x_{kl}^{j\alpha}m_q\ot m_{\alpha}\ot m_j\Bigl)\\
&=&\sum_{j,\alpha, p,i}x_{kl}^{j\alpha}y_{\alpha q}^{ip}
m_p\ot m_{i}\ot m_j\\
&=&\sum_{i,j,p}\Bigl(\sum_{\alpha}x_{kl}^{j\alpha}y_{\alpha q}^{ip} \Bigl)
m_p\ot m_{i}\ot m_j
\end{eqnarray*}
Hence, the conclusion follows.
\end{proof}

In the next proposition we shall evidence a few equations which are
equivalent to the ${\cal D}$-equation.

\begin{proposition}
Let $M$ be a vector space and $R\in \End_k(M\ot M)$. The following
statements are equivalent:
\begin{enumerate}
\item $R$ is a solution of the ${\cal D}$-equation.
\item $T:=R\tau$ is a solution of the equation:
$$T^{12}T^{13}=T^{23}T^{13}\tau^{(123)}$$
\item $U:=\tau R$ is a solution of the equation:
$$U^{13}U^{23}=\tau^{(123)}U^{13}U^{12}$$
\item $W:=\tau R\tau$ is a solution of the equation:
$$\tau^{(123)}W^{23}W^{13}=W^{12}W^{13}\tau^{(123)}$$
\end{enumerate}
\end{proposition}

\begin{proof}
1 $\Leftrightarrow$ 2. As $R=T\tau$, $R$ is a solution of the
${\cal D}$-equation if and only if
\begin{equation}\label{t}
T^{12}\tau^{12}T^{23}\tau^{23}=T^{23}\tau^{23}T^{12}\tau^{12}.
\end{equation}
But
$$
\tau^{12}T^{23}\tau^{23}=T^{13}\tau^{13}\tau^{12},\quad \mbox{and} \quad
\tau^{23}T^{12}\tau^{12}=T^{13}\tau^{12}\tau^{13}.
$$
Hence, the equation (\ref{t}) is equivalent to
$$T^{12}T^{13}\tau^{13}\tau^{12}=T^{23}T^{13}\tau^{12}\tau^{13}.$$
Now, the conclusion follows because
$\tau^{12}\tau^{13}\tau^{12}\tau^{13}=\tau^{(123)}$.

1 $\Leftrightarrow$ 3. $R=\tau U$. Hence, $R$ is a solution of the
${\cal D}$-equation if and only if
\begin{equation}\label{ti}
\tau^{12}U^{12}\tau^{23}U^{23}=\tau^{23}U^{23}\tau^{12}U^{12}.
\end{equation}
But
$$
\tau^{12}U^{12}\tau^{23}=\tau^{23}\tau^{13}U^{13},\quad \mbox{and} \quad
\tau^{23}U^{23}\tau^{12}=\tau^{23}\tau^{12}U^{13}.
$$
Hence, the equation (\ref{ti}) is equivalent to
$$U^{13}U^{23}=\tau^{13}\tau^{12}U^{13}U^{12}$$
and we are done as $\tau^{13}\tau^{12}=\tau^{(123)}$.

1 $\Leftrightarrow$ 4. $R=\tau W\tau$. It follows that $R$ is a solution
of the ${\cal D}$-equation if and only if
\begin{equation}\label{tii}
\tau^{12}W^{12}\tau^{12}\tau^{23}W^{23}\tau^{23}=
\tau^{23}W^{23}\tau^{23}\tau^{12}W^{12}\tau^{12}
\end{equation}
Using the formulas
$$
\tau^{12}\tau^{23}=\tau^{13}\tau^{12}, \quad
\tau^{23}\tau^{12}=\tau^{13}\tau^{23},
$$
$$
\tau^{12}W^{12}\tau^{13}=\tau^{12}\tau^{13}W^{23}, \quad
\tau^{12}W^{23}\tau^{23}=W^{13}\tau^{12}\tau^{23}
$$
$$
\tau^{23}W^{23}\tau^{13}=\tau^{23}\tau^{13}W^{12}, \quad
\tau^{23}W^{12}\tau^{12}=W^{13}\tau^{23}\tau^{12}
$$
the equation (\ref{tii}) is equivalent to
$$
\tau^{12}\tau^{13}W^{23}W^{13}\tau^{12}\tau^{23}=
\tau^{23}\tau^{13}W^{12}W^{13}\tau^{23}\tau^{12}.
$$
The prove is done as
$$
\tau^{13}\tau^{23}\tau^{12}\tau^{13}=\tau^{23}\tau^{12}\tau^{23}\tau^{12}=
\tau^{(123)}.
$$
\end{proof}

\begin{examples}
1. Suppose that $R\in \End_k(M\ot M)$ is bijective. Then, $R$ is a
solution of the ${\cal D}$-equation if and only if $R^{-1}$ is.

2. Let $(m_i)_{i\in I}$ be a basis of $M$ and $(a_{ij})_{i,j\in I}$ be a
family of scalars of $k$. Then, $R:M\ot M\to M\ot M$,
$R(m_i\ot m_j)=a_{ij}m_i\ot m_j$, for all $i$, $j\in I$, is a solution of
the ${\cal D}$-equation. In particular, the identity map $Id_{M\ot M}$ is a
solution of the ${\cal D}$-equation.

3. Let $M$ be a finite dimensional vector space and $u$ an automorphism
of $M$. If $R$ is a solution of the ${\cal D}$-equation then
${}^{u}R:= (u\ot u)R(u\ot u)^{-1}$ is also a solution of the
${\cal D}$-equation.

4. Let $f$, $g\in \End_k(M)$ and $R=f\ot g$. Then $R$ is a solution
of the ${\cal D}$-equation if and only if $fg=gf$.

Indeed, a direct computation shows that
$$
R^{23}R^{12}=f\ot fg\ot g, \quad  R^{12}R^{23}=f\ot gf\ot g
$$
i.e. the conclusion follows.

In particular, we consider the example considered in \cite{LR} (pg. 339).
Let $M$ be a two-dimensional vector space with $\{m_1, m_2 \}$ a basis.
Let $f$, $g\in \End_k(M)$ such that with respect to the given basis are:
$$
f=
\left(
\begin{array}{cc}
a&1\\
0&a
\end{array}
\right), \quad
g=
\left(
\begin{array}{cc}
b&c\\
0&b
\end{array}
\right)
$$
where $a$, $b$, $c$ are scalars of $k$. Then, $R=f\ot g$, with respect
to the ordonate basis $\{m_1\ot m_1, m_1\ot m_2, m_2\ot m_1, m_2\ot m_2 \}$
of $M\ot M$, is given by
\begin{equation}\label{100}
R=
\left(
\begin{array}{cccc}
ab&ac&b&c\\
0&ab&0&b\\
0&0&ab&ac\\
0&0&0&ab
\end{array}
\right)
\end{equation}
Then $R$ is a solution for both the quantum Yang-Baxter equation and the
${\cal D}$-equation.

5. Let $R\in {\cal M}_4(k)$ given by
$$
R=
\left(
\begin{array}{cccc}
a&0&0&0\\
0&b&c&0\\
0&d&e&0\\
0&0&0&f
\end{array}
\right)
$$
It can be proven, by a direct computation, that $R$ is a solution of the
${\cal D}$-equation if and only if $c=d=0$. In particular, if $q\in k$,
$q\neq 0$, $q\neq 1$, the two dimensional Yang-Baxter operator
$$
R=
\left(
\begin{array}{cccc}
q&0&0&0\\
0&1&q-q^{-1}&0\\
0&0&1&0\\
0&0&0&q
\end{array}
\right)
$$
is a solution of the quantum Yang-Baxter equation and is not a solution
for the ${\cal D}$-equation.

6. Let $G$ be a group and $M$ be a left $k[G]$-module.
Suppose that there exists
$\{M_{\sigma }\mid \sigma \in G \}$ a family of $k[G]$-submodules
of $M$ such that $M=\oplus_{\sigma \in G}M_{\sigma}$.
If $m\in M$, then $m$ is a finite sum of homogenous elements
$m=\sum m_{\sigma}$. The map
\begin{equation}\label{gr}
R:M\ot M\to M\ot M, \quad R(n\ot m)=
\sum_{\sigma}\sigma\cdot n\ot m_{\sigma},\; \forall n, m\in M
\end{equation}
is a solution of the ${\cal D}$-equation. Furthermore, if $G$ is non-abelian,
$R$ is not a solution of the quantum Yang-Baxter equation.

Indeed, it is enought to prove that (\ref{Ceq}) holds only for homogenous
elements. Let $m_{\sigma}\in M_{\sigma}$, $m_{\tau}\in M_{\tau}$ and
$m_{\theta}\in M_{\theta}$. Then,
\begin{eqnarray*}
R^{23}R^{12}(m_{\sigma}\ot m_{\tau}\ot m_{\theta})&=&
R^{23}(\tau\cdot m_{\sigma}\ot m_{\tau}\ot m_{\theta})\\
&=&\tau\cdot m_{\sigma}\ot \theta\cdot m_{\tau}\ot m_{\theta}
\end{eqnarray*}
and
\begin{eqnarray*}
R^{12}R^{23}(m_{\sigma}\ot m_{\tau}\ot m_{\theta})&=&
R^{12}(m_{\sigma} \ot \theta \cdot m_{\tau} \ot m_{\theta})\\
\text{($\theta\cdot m_{\tau}\in M_{\tau}$)}
&=&\tau\cdot m_{\sigma}\ot \theta\cdot m_{\tau}\ot m_{\theta}
\end{eqnarray*}
Hence $R$ is a solution of the ${\cal D}$-equation. On the other hand
$$
R^{12}R^{13}R^{23}(m_{\sigma}\ot m_{\tau}\ot m_{\theta})=
\tau\theta\cdot m_{\sigma}\ot \theta\cdot m_{\tau}\ot m_{\theta}
$$
and
$$
R^{23}R^{13}R^{12}(m_{\sigma}\ot m_{\tau}\ot m_{\theta})=
\theta\tau\cdot m_{\sigma}\ot \theta\cdot m_{\tau}\ot m_{\theta}
$$
Hence, if $G$ is a non-abelian group, then $R$ is not a solution of the
quantum Yang-Baxter equation.
\end{examples}

\section{Categories of Long dimodules and the ${\cal D}$-equation}

We shall start this section by adapting a definition given by
Long in \cite{Lo} for a commutative and cocommutative bialgebra to the
case of an arbitrary bialgebra.

\begin{definition}
Let $H$ be a bialgebra. A Long $H$-dimodule is a
threetuple $(M,\cdot ,\rho)$, where $(M,\cdot)$ is a left $H$-module,
$(M,\rho)$ is a right $H$-comodule such that the following
compatibility condition holds:
\begin{equation}\label{C}
\rho(h\cdot m)=\sum h\cdot m_{<0>}\ot m_{<1>}
\end{equation}
for all $h\in H$ and $m\in M$.
\end{definition}
The category of $H$-dimodules over $H$ and $H$-linear $H$-colinear
maps will be denoted by ${}_H{\cal L}^H$.

\begin{examples}
1. Let $G$ be a group and $H=k[G]$ be the groupal Hopf algebra. Let
$(M,\cdot)$ be a left $k[G]$-module. Then, $(M,\cdot ,\rho)$ is an object
in ${}_{k[G]}{\cal L}^{k[G]}$ if and only if there exists
$\{M_{\sigma }\mid \sigma \in G \}$ a family of $k[G]$-submodules
of $M$ such that $M=\oplus_{\sigma \in G}M_{\sigma}$.

Indeed, $(M,\rho)$ is a right $k[G]$-comodule if and only if there exists
$\{M_{\sigma }\mid \sigma \in G \}$ a family of $k$-subspaces
of $M$ such that $M=\oplus_{\sigma \in G}M_{\sigma}$ (see \cite{M}).
Recall that, for $\sigma \in G$,
$$m_{\sigma}\in M_{\sigma}\quad \mbox{if and only
if}\quad \rho(m_{\sigma})=m_{\sigma}\ot \sigma.$$
Now, let $g\in G$ and
$m_{\sigma}\in M_{\sigma}$. The compatibility condition (\ref{C}) turn into
$\rho(g\cdot m_{\sigma})=g\cdot m_{\sigma}\ot \sigma$ which is equivalent to
$g\cdot M_{\sigma}\subseteq M_{\sigma}$ for any $g\in G$, i.e. $M_{\sigma}$
is a left $k[G]$-submodule of $M$ for all $\sigma \in G$.

In fact, we have proven that the categories ${}_{k[G]}{\cal L}^{k[G]}$ and
$\Bigl({}_{k[G]}{\cal M}\Bigl)^{(G)}$ are equivalent.

Let us now suppose that $M$ is a decomposable representation on $G$
with the Long length smaller or equal to the cardinal of $G$. Let $X$ be a
subset of $G$ and $\{M_{x}\mid x \in X \}$ a family of idecomposable
$k[G]$-submodules of $M$ such that $M=\oplus_{x \in X}M_{x}$. Then
$M$ is a right comodule over $k[X]$, and as $X\subseteq G$, $M$ can be
viewed as a right $k[G]$-comodule. Hence, $M\in {}_{k[G]}{\cal L}^{k[G]}$.
We obtain that the category of representations on $G$ with the Long length
smaller or equal to the cardinal of $G$ can be viewed as a subcategory of
${}_{k[G]}{\cal L}^{k[G]}$.

2. Let $(N,\cdot)$ be a left $H$-module. Then $N\ot H$ is an object in
${}_H{\cal L}^H$ with the structures:
$$
h\bullet (n\ot l):=h\cdot n\ot l, \qquad
\rho(n\ot l)=\sum n\ot l_{(1)}\ot l_{(2)}
$$
for all $h$, $l\in H$, $n\in N$. Hence, we have constructed a functor
$$
\bullet\ot H:{}_H{\cal M}\to {}_H{\cal L}^H.
$$
3. Let $(M,\rho)$ be a right $H$-comodule. Then $H\ot M$ is an object in
${}_H{\cal L}^H$ with the structures:
$$
h\bullet (l\ot m):=hl\ot m, \qquad
\rho_{H\ot M}(l\ot m):=\sum l\ot m_{<0>}\ot m_{<1>}
$$
for all $h$, $l\in H$, $m\in M$. Hence, we obtain a functor
$$
H\ot \bullet:{\cal M}^H\to {}_H{\cal L}^H.
$$
4. Let $(N,\cdot)$ be a left $H$-module. Then, with the trivial structure
of right $H$-comodule, $\rho :N\to N\ot H$, $\rho(n):=n\ot 1$, for all
$n\in N$, $(N,\cdot, \rho)$ is an object in ${}_H{\cal L}^H$.

5. Let $(M,\rho)$ be a right $H$-comodule. Then, with the trivial structure
of left $H$-module, $h\cdot m:=\varepsilon(h)m$, for all $h\in H$, $m\in M$,
$(M,\cdot, \rho)$ is an object in ${}_H{\cal L}^H$.
\end{examples}

\begin{remarks}
1. In the study of the category ${}_H{\cal L}^H$ an important role is
played by the forgetful functors:
$$
F^H :{}_H{\cal L}^H \to {}_H{\cal M}\quad \mbox{and} \quad
{}_HF :{}_H{\cal L}^H \to {\cal M}^H
$$
where $F^H$ (resp. ${}_HF$) is the functor forgetting the $H$-comodule
(resp. $H$-module) structure.

The principle according to which forgetful type functors have adjoints is
also valid for the category ${}_H{\cal L}^H$. We have:
\begin{enumerate}
\item The functor $H\ot \bullet:{\cal M}^H\to {}_H{\cal L}^H$ is a left
adjoint of the functor ${}_HF :{}_H{\cal L}^H \to {\cal M}^H$.
\item The functor $\bullet\ot H:{}_H{\cal M}\to {}_H{\cal L}^H$ is a right
adjoint of the functor $F^H :{}_H{\cal L}^H \to {}_H{\cal M}$.
\end{enumerate}
The adjoint pair $(H\ot \bullet,\; {}_HF)$ is given by the following natural
transformation:
$$
\eta :Id_{{\cal M}^H}\to {}_HF\circ(H\ot \bullet),\;
\eta_N:N\to H\ot N,\; \eta_N(n):=1\ot n,\; \forall n\in N
$$
$$
\theta :(H\ot \bullet)\circ {}_HF\to Id_{{}_H{\cal L}^H},\:
\theta_M:H\ot M\to M,\: \theta_M(h\ot m):=hm,\: \forall h\in H, m\in M
$$
where $N\in {\cal M}^H$, $M\in {}_H{\cal L}^H$.

The adjoint pair $(F^H, \bullet \ot H)$ is given by the following natural
transformation:
$$
\eta :Id_{{}_H{\cal L}^H}\to (\bullet\ot H)\circ F^H,\:
\eta_M:M\to M\ot H,\: \eta_M(m):=\sum m_{<0>}\ot m_{<1>},\: \forall m\in M
$$
$$
\theta :F^H\circ(\bullet\ot H)\to Id_{{}_H{\cal M}},\;
\theta_N:N\ot H\to N,\; \theta_N(n\ot h):=\varepsilon(h)n,\;
\forall h\in H, n\in N
$$
where  $M\in {}_H{\cal L}^H$, $N\in {}_H{\cal M}$.

Now, using the general properties of pairs of adjoint functors we obtain
(see Corollary 2.7 of \cite{CMZ1} for a similar results regarding the
Yetter-Drinfel'd category ${}_H{\cal YD}^H$):
\begin{enumerate}
\item the functor $\bullet\otimes H: {}_H{\cal M}\to {}_H{\cal L}^H$
preserves the injective cogenerators. In particular,
$\Hom_{\nint}(H,\nrat /\nint)\otimes H$ is an injective cogenerator
of ${}_ H{\cal L}^H$.
\item the functor $H\otimes\bullet: {\cal M}^H\to {}_H{\cal L}^H$ preserves
generators. In particular, if $H$ is cofrobenius as a coalgebra, then
$H\otimes H$ is a projective generator of ${}_H{\cal L}^H$.
\end{enumerate}

2. If $M$, $N\in {}_H{\cal L}^H$, then $M\ot N$ is also an object in
${}_H{\cal L}^H$ with the natural structures
$$
h\bullet (m\ot n)=\sum h_{(1)}\cdot m\ot h_{(2)}\cdot n,\quad
\rho(m\ot n)=\sum m_{<0>}\ot n_{<0>}\ot m_{<1>}n_{<1>}
$$
for all $h\in H$, $m\in M$, $n\in N$. $k$ can be view as an object in
${}_H{\cal L}^H$ with the trivial structures
$$
h\cdot a=\varepsilon(h)a,\qquad \rho(a)=a\ot 1
$$
for all $h\in H$, $a\in k$. It is easy to see that
$({}_H{\cal L}^H, \ot , k)$ is a monoidal category. This can also be
obtained from Proposition 2.1 of \cite{COZ}, keeping in mind that the
category ${}_H{\cal L}^H$ is a special case of the Doi-Hopf modules
category.

3. If $M\in {}_H{\cal L}^H$, then $M$ is a left $H\ot H^*$-module with the
structure
$$
(h\ot h^*)\cdot m:=\sum <h^*, m_{<1>}> hm_{<0>}
$$
for all $h\in H$, $h^*\in H^*$ and $m\in M$. Furthermore, if $H$ is
finite dimensional the categories ${}_H{\cal L}^H$ and
${}_{H\ot H^*}{\cal M}$ are equivalent.

4. The category ${}_H{\cal L}^H$ can be generalized to a category
${}_A{\cal L}^C$.
Let $A$ be an algebra and $C$ be a coalgebra over $k$. Then an object in
${}_A{\cal L}^C$ is a threetuple $(M,\cdot ,\rho)$, where $(M,\cdot)$
is a left $A$-module, $(M,\rho)$ is a right $C$-comodule such that the
following compatibility condition holds:
\begin{equation}\label{Cgen}
\rho(a\cdot m)=\sum a\cdot m_{<0>}\ot m_{<1>}
\end{equation}
for all $a\in A$ and $m\in M$.

Let $X$ be a set and $C=k[X]$ be the groupal coalgebra. Then
$(M,\cdot ,\rho)$ is an object
in ${}_A{\cal L}^{k[X]}$ if and only if there exists
$\{M_{x}\mid x \in X \}$ a family of $A$-submodules
of $M$ such that $M=\oplus_{x \in X}M_{x}$.
Hence, the category ${}_A{\cal L}^C$ can be seen as a generalization of
the category of decomposable left $A$-modules.
\end{remarks}

The next proposition evidences the role which can be played by the
category ${}_H{\cal L}^H$ in solving the ${\cal D}$-equation.

\begin{proposition}
Let $H$ be a bialgebra and $(M,\cdot, \rho)$ be a Long $H$-dimodule.
Then the natural map
$$R_{(M,\cdot ,\rho )}(m\ot n)=\sum n_{<1>}\cdot m\ot n_{<0>}$$
is a solution of the ${\cal D}$-equation.
\end{proposition}

\begin{proof}
Let $R=R_{(M,\cdot ,\rho )}$. For $l$, $m$, $n\in M$ we have
\begin{eqnarray*}
R^{12}R^{23}(l\ot m\ot n)&=&
R^{12}\Bigl ( \sum l\ot n_{<1>}\cdot m \ot n_{<0>} \Bigl )\\
&=&\sum (n_{<1>}\cdot m)_{<1>}\cdot l \ot (n_{<1>}\cdot m)_{<0>}\ot n_{<0>}\\
\text{(\mbox{using} (\ref{C}))}
&=&\sum m_{<1>}\cdot l\ot n_{<1>}\cdot m_{<0>}\ot n_{<0>}\\
&=&R^{23}\Bigl(\sum m_{<1>}\cdot l \ot m_{<0>}\ot n \Bigl)\\
&=&R^{23}R^{12}(l\ot m\ot n)
\end{eqnarray*}
i.e. $R$ is a solution of the ${\cal D}$-equation.
\end{proof}

\begin{lemma} \label{trei}
Let $H$ be a bialgebra, $(M,\cdot)$ a left $H$-module
and $(M,\rho)$ a right $H$-comodule. Then the set
$$
\{h\in H\mid \rho(h\cdot m)=\sum h\cdot m_{<0>}\ot m_{<1>},
\forall m\in M \}
$$
is a subalgebra of $H$.
\end{lemma}

\begin{proof} Straightforward.
\end{proof}

We obtain from this lemma that if a left $H$-module and  right
$H$-comodule $M$ satisfies the condition of compatibility (\ref{C})
for a set of generators as an algebra of $H$ and for a basis of $M$, then
$M$ is a Long $H$-dimodule.

Now we shall prove the main results of this section. By a FRT type
theorem, we shall prove that in the finite dimensional case, any solution of
the ${\cal D}$-equation has the form $R_{(M,\cdot,\rho)}$, where $(M,\cdot,\rho)$
is a Long $D(R)$-dimodule over a special bialgebra $D(R)$.

\begin{theorem}
Let $M$ be a finite dimensional vector space and $R\in \End_k(M\ot M)$
be a solution of the ${\cal D}$-equation. Then:
\begin{enumerate}
\item There exists a bialgebra $D(R)$ such that $(M,\cdot,\rho)$ has a
structure of object in ${}_{D(R)}{\cal L}^{D(R)}$  and
$R=R_{(M,\cdot, \rho)}$.
\item The bialgebra $D(R)$ is a universal object with this property:
if $H$ is a bialgebra such that
$(M,\cdot^{\prime}, \rho^{\prime})\in {}_H{\cal L}^H$ and
$R=R_{(M,\cdot^{\prime}, \rho^{\prime})}$ then there exists a unique
bialgebra map $f:D(R)\to H$ such that $\rho^{\prime}=(I\ot f)\rho$.
Furthermore, $a\cdot m=f(a)\cdot^{\prime}m$, for all $a\in D(R)$,
$m\in M$.
\end{enumerate}
\end{theorem}

\begin{proof}
$1$. Let $\{m_1,\cdots ,m_n \}$
be a basis for $M$ and $(x_{uv}^{ji})_{i,j,u,v}$ a family
of scalars of $k$ such that
\begin{equation}
R(m_v\ot m_u)=\sum_{i,j}x_{uv}^{ji}m_i\ot m_j
\end{equation}
for all $u$, $v=1,\cdots ,n$.

Let $(C, \Delta, \varepsilon)={\cal M}^n(k)$, be the comatrix coalgebra
of order $n$.  Let $\rho :M\to M\ot C$ given by
\begin{equation}\label{ro}
\rho(m_l)=\sum_{v=1}^{n}m_v\ot c_{vl}
\end{equation}
for all $l=1,\cdots, n$. Then, $M$ is a right $C$-comodule.
Let $T(C)$ be the unique bialgebra structure on the tensor
algebra $T(C)$ which extends $\Delta$ and $\varepsilon$.
As the inclusion $i:C\to T(C)$ is a coalgebra map, $M$ has a right
$T(C)$-comodule structure via
$$
M\stackrel{\rho}{\longrightarrow}
M\ot C\stackrel{I\ot i}{\longrightarrow}M\ot T(C)
$$
There will be no confusion if we also denote the right
$T(C)$-comodule structure on $M$ with $\rho$.

Now, we will put a left $T(C)$-module structure on $M$ in such a way that
$R=R_{(M,\cdot, \rho)}$. First we define the $k$-bilinear map
$$
\mu :C\ot M\to M, \quad
\mu(c_{ju}\ot m_v):=\sum_{i} x_{uv}^{ji}m_i
$$
for all $j$, $u$, $v=1,\cdots n$. Using the universal property of the tensor
algebra $T(C)$, there exists a unique left $T(C)$-module structure on
$(M,\cdot )$ such that
$$
c_{ju}\cdot m_v=\mu(c_{ju}\ot m_v)=\sum_{i} x_{uv}^{ji}m_i
$$
for all $j$, $u$, $v=1,\cdots ,n$. For $m_v$, $m_u$ the elements of
the given basis, we have:
\begin{eqnarray*}
R_{(M,\cdot, \rho)}(m_v\ot m_u)&=&
\sum_j c_{ju}\cdot m_v\ot m_j\\
&=&\sum_{i,j} x_{uv}^{ji}m_i\ot m_j\\
&=&R(m_v\ot m_u)
\end{eqnarray*}
Hence, $(M,\cdot, \rho)$ has a structure of left $T(C)$-module and
right $T(C)$-comodule such that $R=R_{(M,\cdot, \rho)}$.

Now, we define the {\sl obstructions} $o(i,j,k,l)$ which measure
how far away $M$ is from a Long $T(C)$-dimodule.
Keeping in mind that $T(C)$ is generated as an algebra by $(c_{ij})$
and using lemma \ref{trei} we compute
$$
\sum h\cdot m_{<0>}\ot m_{<1>} - \rho (h\cdot m)
$$
only for $h=c_{jk}$, and $m=m_l$, for $j$, $k$, $l=1,\cdots, n$. We have:
\begin{eqnarray*}
\sum h\cdot m_{<0>}\ot m_{<1>}&=&
\sum c_{jk}\cdot (m_{l})_{<0>}\ot (m_{l})_{<1>}\\
&=&\sum_{v} c_{jk}\cdot m_v \ot c_{vl}\\
&=&\sum_{v,i}x_{kv}^{ji}m_i\ot c_{vl}\\
&=&\sum_{i}m_i\ot \Bigl(\sum_v x_{kv}^{ji}c_{vl}\Bigl)
\end{eqnarray*}
and
\begin{eqnarray*}
\rho(h\cdot m)&=&\rho (c_{jk}\cdot m_l)\\
&=&\sum_{\alpha}x_{kl}^{j\alpha}(m_{\alpha})_{<0>}\ot (m_{\alpha})_{<1>}\\
&=&\sum_{i,\alpha}x_{kl}^{j\alpha}m_i\ot c_{i\alpha}\\
&=&\sum_im_i\ot \Bigl(\sum_{\alpha}x_{kl}^{j\alpha}c_{i\alpha}\Bigl)
\end{eqnarray*}
Let
\begin{equation}\label{obst}
o(i,j,k,l):=\sum_{v}x_{kv}^{ji}c_{vl} -
\sum_{\alpha}x_{kl}^{j\alpha}c_{i\alpha}
\end{equation}
for all $i$, $j$, $k$, $l=1,\cdots, n$. Then
\begin{equation}\label{pais}
\sum h\cdot m_{<0>}\ot m_{<1>}-
\rho(h\cdot m)=\sum_i m_i\ot o(i,j,k,l)
\end{equation}
Let $I$ be the two-sided ideal of $T(C)$ generated by all $o(i,j,k,l)$,
$i$, $j$, $k$, $l=1,\cdots, n$. Then:

{\sl $I$ is a bi-ideal of $T(C)$ and $I\cdot M=0$.}

We first prove that $I$ is also a coideal and this will result from the
following formula:
\begin{equation}\label{prima}
\Delta(o(i,j,k,l))=\sum_{u}\Bigl( o(i,j,k,u)\ot c_{ul}+
c_{iu}\ot o(u,j,k,l)\Bigl)
\end{equation}
First, we denote that the formula (\ref{obst}) can be written
$$
o(i,j,k,l):=\sum_{v}\Bigl( x_{kv}^{ji}c_{vl}-x_{kl}^{jv}c_{iv} \Bigl)
$$
We have:
\begin{eqnarray*}
\Delta(o(i,j,k,l))&=&\sum_{u,v}
\Bigl( x_{kv}^{ji}c_{vu}\ot c_{ul}-
x_{kl}^{jv}c_{iu}\ot c_{uv}\Bigl)\\
&=&\sum_{u}\Bigl(\sum_v x_{kv}^{ji}c_{vu}\Bigl)\ot c_{ul}-
\sum_{u}c_{iu}\ot\Bigl(\sum_vx_{kl}^{jv} c_{uv}\Bigl)\\
&=&\sum_{u}\Bigl(o(i,j,k,u)+
\sum_{v}x_{ku}^{jv}c_{iv}\Bigl)\ot c_{ul}\\
&-&\sum_{u}c_{iu}\ot \Bigl(-o(u,j,k,l)+
\sum_{v}x_{kv}^{ju}c_{vl} \Bigl) \\
&=&\sum_{u}\Bigl( o(i,j,k,u)\ot c_{ul}+
c_{iu}\ot u(u,j,k,l)\Bigl)
\end{eqnarray*}
where in the last equality we use the fact that
$$
\sum_{u,v}x_{ku}^{jv}c_{iv}\ot c_{ul}=
\sum_{u,v}x_{kv}^{ju}c_{iu}\ot c_{vl}
$$
Hence, the formula (\ref{prima}) holds. On the other hand
$$
\varepsilon\Bigl(o(i,j,k,l)\Bigl)=x_{kl}^{ji}-x_{kl}^{ji}=0
$$
so we proved that $I$ is a coideal of $T(C)$.

Now, we shall prove that $I\cdot M=0$ using the fact that
$R$ is a solution of the ${\cal D}$-equation. Let $n\in M$ and
$j$, $k=1,\cdots ,n$. We have
\begin{equation}\label{adoua}
\Bigl(R^{23}R^{12}-R^{12}R^{23}\Bigl)(n\ot m_k\ot m_j)=
\sum_{r,s}o(r,s,j,k)\cdot n\ot m_r\ot m_s
\end{equation}
Let us compute
\begin{eqnarray*}
\Bigl(R^{23}R^{12} \Bigl)(n\ot m_k\ot m_j)&=&
R^{23}(\sum_{\alpha}c_{\alpha k}\cdot n\ot m_{\alpha}\ot m_{j} )\\
&=&\sum_{\alpha,s}
c_{\alpha k}\cdot n\ot c_{sj}\cdot m_{\alpha}\ot m_{s} \\
&=&\sum_{\alpha,r,s}c_{\alpha k}\cdot n\ot x_{j\alpha}^{sr}m_{r}\ot m_{s}\\
&=&\sum_{r,s}\Bigl(\sum_{\alpha}
x_{j\alpha}^{sr}c_{\alpha k}\Bigl )\cdot n\ot m_{r}\ot m_{s}
\end{eqnarray*}
On the other hand
\begin{eqnarray*}
\Bigl(R^{12}R^{23}\Bigl)(n\ot m_k\ot m_j)&=&
R^{12}(\sum_{s}n\ot c_{sj}\cdot m_k \ot m_s)\\
&=&R^{12}(\sum_{s,\alpha}n\ot x_{jk}^{s\alpha}m_{\alpha}\ot m_s)\\
&=&\sum_{r,s,\alpha}x_{jk}^{s\alpha}c_{r\alpha}\cdot n\ot m_r\ot m_s
\end{eqnarray*}
It follows that
\begin{eqnarray*}
\Bigl(R^{23}R^{12}-R^{12}R^{23}\Bigl)(n\ot m_k\ot m_j)
&=&\sum_{r,s}\Bigl(\sum_{\alpha}x_{j\alpha}^{sr}c_{\alpha k}-
\sum_{\alpha}x_{jk}^{s\alpha}c_{r\alpha} \Bigl)\cdot n\ot m_r\ot m_s\\
&=&\sum_{r,s}o(r,s,j,k)\cdot n\ot m_r\ot m_s
\end{eqnarray*}
i.e. the formula (\ref{adoua}) holds. But $R$ is a solution of the
${\cal D}$-equation, hence $o(r,s,j,k)\cdot n =0$, for all $n\in M$,
$j$, $k$, $r$, $s=1,\cdots, n$. We conclude that
$I$ is a bi-ideal of $T(C)$ and $I\cdot M=0$. Define now
$$
D(R)=T(C)/I
$$
$M$ has a right $D(R)$-comodule structure via the canonical projection
$T(C)\to D(R)$ and a left $D(R)$-module structure as $I\cdot M=0$.
As $(c_{ij})$ generate $D(R)$ and in $D(R)$, $o(i,j,k,l)=0$, for all
$i$, $j$, $k$, $l=1,\cdots, n$, using (\ref{pais}) we get that
$(M,\cdot, \rho)\in {}_{D(R)}{\cal L}^{D(R)}$ and $R=R_{(M,\cdot,\rho)}$.

$2$. Let $H$ be a bialgebra and suppose that
$(M,\cdot^{\prime}, \rho^{\prime})\in {}_H{\cal L}^H$ and
$R=R_{(M,\cdot^{\prime}, \rho^{\prime})}$.
Let $(c_{ij}^{\prime})_{i,j=1,\cdots ,n}$ be a family of elements of $H$
such that
$$
\rho^{\prime}(m_l)=\sum_v m_v\ot c_{vl}^{\prime}
$$
Then
$$
R(m_v\ot m_u)=\sum_j c_{ju}^{\prime}\cdot^{\prime}m_v \ot m_j
$$
and
$$
c_{ju}^{\prime}\cdot^{\prime}m_v=\sum_i x_{uv}^{ji}m_i=c_{ju}\cdot m_v.
$$
Let
$$
o^{\prime}(i,j,k,l)=
\sum_{v}x_{kv}^{ji}c_{vl}^{\prime} -
\sum_{\alpha}x_{kl}^{j\alpha}c_{i\alpha}^{\prime}
$$
From the universal property of the tensor algebra $T(C)$, there exists
a unique algebra map $f_1: T(C)\to H$ such that
$f_1(c_{ij})=c_{ij}^{\prime}$, for all $i$, $j=1,\cdots, n$.
As $(M,\cdot^{\prime}, \rho^{\prime})\in {}_H{\cal L}^H$ we get that
$o^{\prime}(i,j,k,l)=0$, and hence $f_1(o(i,j,k,l))=0$,
for all $i$, $j$, $k$, $l=1,\cdots, n$. So the map $f_1$ factorizes to
the map
$$
f:D(R)\to H, \quad f(c_{ij})=c_{ij}^{\prime}
$$
Of course, for $m_l$ an arbitrary element of the given basis of $M$, we have
$$
(I\ot f)\rho(m_l)=\sum_v m_v\ot f(c_{vl})=\sum_v m_v\ot c_{vl}^{\prime}=
\rho^{\prime}(m_l)
$$
Conversely, the relation $(I\ot f)\rho=\rho^{\prime}$ necessarily
implies $f(c_{ij})=c_{ij}^{\prime}$, which proves the uniqueness of $f$.
This completes the proof of the theorem.
\end{proof}

\begin{remark}
The obstruction elements $o(i,j,k,l)$ are different from the
homogenous elements $d(i,j,k,l)$ defined in \cite{R1} which correspond
to the quantum Yang-Baxter equation, and are also different from the
obstruction elements $\chi(i,j,k,l)$ which appear in \cite{M1} in connexion
with the Hopf equation. The main difference consists in the fact that in the
graded algebra $T({\cal M}^n(k))$ the elements $o(i,j,k,l)$
are of degree one, i.e. are elements of the comatrix coalgebra
${\cal M}^n(k)$. This will lead us in the next section to the study of
some special functions defined only for a coalgebra, which will also play an
important role in solving the ${\cal D}$-equation.
\end{remark}

\begin{examples}
1. Let $a$, $b$, $c\in k$ and $R\in {\cal M}_4(k)$ given by equation
(\ref{100}), which is a solution for both the quantum Yang-Baxter equation
and the ${\cal D}$-equation. In \cite{LR}, if $ab\neq 0$, $ac+b\neq 0$ a
labour-intensive computation will give a description of the bialgebra
$A(R)$, obtained looking at $R$ as a solution of the quantum
Yang-Baxter equation.

Below we shall describe the bialgebra $D(R)$, which is obtained considering
$R$ as a solution for the ${\cal D}$-equation.
If $(b,c)=(0,0)$ then $R=0$ i.e. $D(R)=T({\cal M}^4(k))$. Suppose now that
$(b,c)\neq (0,0)$. If we write
$$
R(m_v\ot m_u)=\sum_{i,j=1}^{2}x_{uv}^{ji}m_i\ot m_j
$$
we get that among the elements $(x_{uv}^{ji})$, the only nonzero elements
are:
\begin{eqnarray}\label{scalari}
x_{11}^{11}=ab,\quad x_{21}^{11}=ac,\quad x_{21}^{21}=ab,
\quad x_{12}^{11}=b, \\
x_{12}^{12}=ab,\quad x_{21}^{11}=c,\quad x_{22}^{21}=b,\quad
x_{22}^{12}=ac,\quad x_{22}^{22}=ab \nonumber
\end{eqnarray}
The sixteen relation $o(i,j,k,l)=0$, written in the lexicografical order
acording to $(i,j,k,l)$, starting with $(1,1,1,1)$ are
$$
abc_{11}+bc_{21}=abc_{11},\quad abc_{12}+bc_{22}=bc_{11}+abc_{12}
$$
$$
acc_{11}+cc_{21}=acc_{11}, \quad acc_{12}+cc_{22}=cc_{11}+acc_{12}
$$
$$
0=0,\quad 0=0,\quad abc_{11}+bc_{21}=abc_{11},\quad
abc_{12}+bc_{22}=bc_{11}+abc_{12}
$$
$$
abc_{21}=abc_{21},\quad abc_{22}=abc_{22},\quad acc_{21}=acc_{21},\quad
acc_{22}=cc_{21}+acc_{22}
$$
$$
0=0,\quad 0=0,\quad abc_{21}=abc_{21}, \quad abc_{22}=abc_{22}+bc_{21}
$$
It remains only the four relation
$$
bc_{21}=0,\quad bc_{22}=bc_{11},\quad
cc_{21}=0,\quad cc_{22}=cc_{11}
$$
As, $(b,c)\neq (0,0)$, there are only two linear independent relation:
$$
c_{21}=0,\qquad c_{22}=c_{11}
$$
Now, if we denote $c_{11}=x$, $c_{12}=y$ we obtain that $D(R)$ can be
described as follows:

$\bullet$ as an algebra $D(R)=k<x, y>$, the free algebra generated
by $x$ and $y$.

$\bullet$ The comultiplication $\Delta$ and the counity $\varepsilon$
are given by
$$
\Delta(x)=x\ot x, \quad \Delta(y)=x\ot y+y\ot x,\quad
\varepsilon(x)=1, \quad \varepsilon(y)=0.
$$
We observe that the bialgebra $D(R)$ does not depend on the parameters
$a$, $b$, $c$.

2. Let $q\in k$ be a scalar and $R_{q}\in {\cal M}_4(k)$ given by
$$
R_{q}=
\left(
\begin{array}{cccc}
0&-q&0&-q^2\\
0&1&0&q\\
0&0&0&0\\
0&0&0&0
\end{array}
\right)
$$
Then $R_{q}$ is a solution of the Hopf equation (cf. \cite{M1}), and of
the ${\cal D}$-equation, as $R_{q}$ has the form $f\ot g$ with $fg=gf$. In
\cite{M1} we have described the bialgebra $B(R_q)$ which arises thinking
of $R_q$ as a solution of the Hopf equation.
The bialgebra $D(R_q)$ obtained by viewing $R_q$ as a solution of the
${\cal D}$-equation has a much simpler description and is independent of $q$.

$\bullet$ as an algebra $D(R_q)=k<x, y>$, the free algebra generated
by $x$ and $y$.

$\bullet$ The comultiplication $\Delta$ and the counity $\varepsilon$
are given in such a way $x$ and $y$ are groupal elements
$$
\Delta(x)=x\ot x, \quad \Delta(y)=y\ot y\quad
\varepsilon(x)=\varepsilon(y)=1.
$$
The above description can be obtained through a computation similar to the
one from the previous example. Among the sixteen relations $o(i,j,k,l)=0$,
the only linear independent ones are:
$$
c_{21}=0, \qquad c_{12}=q(c_{11}-c_{22})
$$
If we denote $c_{11}=x$, $c_{22}=y$ the conclusion follows.

3. We shall now present an example which has a geometric flavour.
Let $f\in \End_k(k^2)$, $f((x,y))=(x,0)$ for all $(x,y)\in k^2$, i.e.
$f$ is the projection of the plane $k^2$ on the $Ox$ axis. With respect to
the canonical basis $\{e_1, e_2\}$ of $k^2$, $f$ has the form
$$
f=
\left(
\begin{array}{cc}
1&0\\
0&0
\end{array}
\right)
$$
$g\in \End_k(k^2)$ commute with $f$ if and only if with respect to
the canonical basis, $g$ has the form
$$
g=
\left(
\begin{array}{cc}
a&0\\
0&b
\end{array}
\right)
$$
where $a$, $b\in k$. Then $R=f\ot g$ with respect to the ordonate basis
$\{e_1\ot e_1, e_1\ot e_2, e_2\ot e_1, e_2\ot e_2\}$ is given by
$$
R=
\left(
\begin{array}{cccc}
a&0&0&0\\
0&b&0&0\\
0&0&0&0\\
0&0&0&0
\end{array}
\right)
$$
We shall describe the bialgebra $D(R)$. Suppose $(a,b)\neq 0$ (otherwise
$R=0$ and $D(R)=T({\cal M}^4(k))$ ).
Among the sixteen relations $o(i,j,k,l)=0$, the only linear independent
ones are:
$$
ac_{12}=ac_{21}=bc_{12}=bc_{21}=0
$$
As $(a,b)\neq 0$ we obtain $c_{12}=c_{21}=0$. Now, if we denote
$c_{11}=x$, $c_{22}=y$ we get that the bialgebra $D(R)$ is exactly the
bialgebra described in the previous example.
\end{examples}

\section{Special functions on a coalgebra and the ${\cal D}$-equation}
In this section we shall introduce the concept of ${\cal D}$-map on a
coalgebra: if $C$ is a coalgebra and $I$ a coideal of $C$, then a
${\cal D}$-map is a k-bilinear map
$\sigma :C\ot C/I\to k$ which satisfies the equation (\ref{spe})
presented below. This condition
ensures that, for any right $C$-comodule $M$, the natural map
$R_{\sigma}$ is a solution for the ${\cal D}$-equation. Conversely, in the
finite dimensional case, any solution for the ${\cal D}$-equation arises in
this way.

If $C$ is a coalgebra and $I$ is a coideal of $C$ then the elements of
the quotient $C/I$ will be denoted by $\overline{c}$. If $(M,\rho)$ is
a right $C$-comodule, then $(M,\overline{\rho})$ is a right $C/I$-comodule
via $\overline{\rho}(m)=\sum m_{<0>}\ot \overline{m_{<1>}}$,
for all $m\in M$.

\begin{definition}
Let $C$ be a coalgebra and $I$ be a coideal of $C$. A $k$-bilinear map
$\sigma :C\ot C/I\to k$ is called a ${\cal D}$-map if
\begin{equation}\label{spe}
\sum \sigma(c_{(1)}\ot \overline{d})\:\overline{c_{(2)}}=
\sum \sigma(c_{(2)}\ot \overline{d})\:\overline{c_{(1)}}
\end{equation}
for all $c\in C$, $\overline{d}\in C/I$. If $I=0$, then $\sigma$ is called
strongly ${\cal D}$-map.
\end{definition}

\begin{examples}
1. If $C$ is cocommutavive then any $k$-bilinear map
$\sigma :C\ot C/I\to k$ is a ${\cal D}$-map.

2. Let $C$ be a coalgebra, $I$ be a coideal of $C$ and
$f\in \Hom_k(C/I, k)$. Then
$$\sigma_f :C\ot C/I\to k,\quad
\sigma_f(c\ot \overline{d}):=\varepsilon(c)f(\overline{d})$$
for all $c\in C$, $\overline{d}\in C/I$, is a ${\cal D}$-map.
In particular, $\sigma :C\ot C/I\to k$,
$\sigma(c\ot \overline{d}):=\varepsilon(c)\varepsilon(\overline{d})$,
for all $c\in C$, $\overline{d}\in C/I$, is a ${\cal D}$-map.

3. Let $C:={\cal M}^n(k)$ be the comatrix coalgebra of order $n$.
Let $a\in k$ be a scalar of $k$. Then the map
$$\sigma :C\ot C\to k,\quad
\sigma (c_{ij}\ot c_{pq}):=\delta_{ij}a$$
for all $i$, $j$, $p$, $q=1,\cdots, n$ is a strongly ${\cal D}$-map.

Indeed, for $c=c_{ij}$, $d=c_{pq}$ we have
$$
\sum \sigma(c_{(1)}\ot d)c_{(2)}=
\sum_{t}\sigma(c_{it}\ot c_{pq})c_{tj}=ac_{ij},
$$
and
$$
\sum \sigma(c_{(2)}\ot d)c_{(1)}=
\sum_{t}\sigma(c_{tj}\ot c_{pq})c_{it}=ac_{ij},
$$
for all $i$, $j$, $p$, $q=1,\cdots, n$.
\end{examples}

\begin{proposition}
Let $C$ be a coalgebra, $I$ be a coideal of $C$ and $\sigma :C\ot C/I\to k$
a ${\cal D}$-map. Let $(M,\rho)$ be a right $C$-comodule. Then, the special
map
$$
R_{\sigma}:M\ot M\to M\ot M, \quad
R_{\sigma}(m\ot n)=
\sum \sigma(m_{<1>}\ot \overline{n_{<1>}})m_{<0>}\ot n_{<0>}
$$
is a solution of the ${\cal D}$-equation.
\end{proposition}

\begin{proof}
Let $l$, $m$, $n\in M$ and put $R=R_{\sigma}$. We have
\begin{eqnarray*}
R^{12}R^{23}(l\ot m\ot n)&=&
R^{12}\Bigl(\sum \sigma(m_{<1>}\ot \overline{n_{<1>}})
l\ot m_{<0>}\ot n_{<0>}\Bigl)\\
&=&\sum \sigma(m_{<1>}\ot \overline{n_{<1>}})
\sigma(l_{<1>}\ot \overline{m_{<0><1>}})l_{<0>}\ot m_{<0><0>}\ot n_{<0>}\\
&=&\sum \sigma(m_{<1>(2)}\ot \overline{n_{<1>}})
\sigma(l_{<1>}\ot \overline{m_{<1>(1)}})l_{<0>}\ot m_{<0>}\ot n_{<0>}\\
&=&\sum \sigma \Bigl(l_{<1>}\ot
\sigma(m_{<1>(2)}\ot \overline{n_{<1>}})\overline{m_{<1>(1)}} \Bigl)
l_{<0>}\ot m_{<0>}\ot n_{<0>}\\
\text{(\mbox{using (\ref{spe})})}
&=&\sum \sigma \Bigl(l_{<1>}\ot
\sigma(m_{<1>(1)}\ot \overline{n_{<1>}})\overline{m_{<1>(2)}} \Bigl)
l_{<0>}\ot m_{<0>}\ot n_{<0>}\\
&=&\sum \sigma(l_{<1>}\ot \overline{m_{<2>}})
\sigma(m_{<1>}\ot \overline{n_{<1>}})l_{<0>}\ot m_{<0>}\ot n_{<0>}
\end{eqnarray*}
and
\begin{eqnarray*}
R^{23}R^{12}(l\ot m\ot n)&=&
R^{23}\Bigl(\sum \sigma(l_{<1>}\ot \overline{m_{<1>}})
l_{<0>}\ot m_{<0>}\ot n\Bigl)\\
&=&\sum \sigma(l_{<1>}\ot \overline{m_{<1>}})
\sigma(m_{<0><1>}\ot \overline{n_{<1>}})l_{<0>}\ot m_{<0><0>}\ot n_{<0>}\\
&=&\sum \sigma(l_{<1>}\ot \overline{m_{<2>}})
\sigma(m_{<1>}\ot \overline{n_{<1>}})l_{<0>}\ot m_{<0>}\ot n_{<0>}
\end{eqnarray*}
i.e. $R_{\sigma}$ is a solution of the ${\cal D}$-equation.
\end{proof}

\begin{theorem}
Let $n$ be a positive integer number, $M$ be a $n$ dimensional vector space
and $R\in \End_k(M\ot M)$ a solution of the ${\cal D}$-equation. Then:
\begin{enumerate}
\item There exist a coideal $I(R)$ of the comatrix coalgebra
${\cal M}^n(k)$ and a unique ${\cal D}$-map
$$\sigma :{\cal M}^n(k)\ot {\cal M}^n(k)/I(R)\to k$$
such that $R=R_{\sigma}$.
Furthermore, if $R$ is bijective, $\sigma$ is invertible in the convolution
algebra $\Hom_k({\cal M}^n(k)\ot {\cal M}^n(k)/I(R),\; k)$.
\item If $R\tau=\tau R$, then there exists a coalgebra $C(R)$ and a unique
strongly ${\cal D}$-map
$$\tilde{\sigma} :C(R)\ot C(R)\to k$$
such that $M$ has a structure of right $C(R)$-comodule and
$R=R_{\tilde{\sigma}}$.
\end{enumerate}
\end{theorem}

\begin{proof}
1. Let $\{m_1,\cdots ,m_n \}$ be a basis for $M$ and
$(x_{uv}^{ji})_{i,j,u,v}$ a family of scalars of $k$ such that
\begin{equation}
R(m_v\ot m_u)=\sum_{i,j}x_{uv}^{ji}m_i\ot m_j
\end{equation}
for all $u$, $v=1,\cdots ,n$. Let ${\cal M}^n(k)$ be the comatrix coalgebra
of order $n$. $M$ has a right ${\cal M}^n(k)$-comodule structure given by
$$\rho(m_l)=\sum_{v=1}^{n}m_v\ot c_{vl}$$
for all $l=1,\cdots, n$. Let $I(R)$ be the $k$-subspace of ${\cal M}^n(k)$
generated by all $o(i,j,k,l)$, $i$, $j$, $k$, $l=1,\cdots,n$. From
equation (\ref{prima}), $I(R)$ is a coideal of ${\cal M}^n(k)$.

First we shall prove the uniqueness.
Let $\sigma :{\cal M}^n(k)\ot {\cal M}^n(k)/I(R)\to k$ be a
${\cal D}$-map such that  $R=R_{\sigma}$. Let $u$, $v=1,\cdots,n$. Then
\begin{eqnarray*}
R_{\sigma}(m_v\ot m_u)&=&\sum \sigma \Bigl((m_v)_{<1>}\ot
\overline{(m_u)_{<1>}}\Bigl)(m_v)_{<0>}\ot (m_u)_{<0>}\\
&=&\sum_{i,j} \sigma(c_{iv}\ot \overline{c_{ju}})m_i\ot m_j
\end{eqnarray*}
Hence $R_{\sigma}(m_v\ot m_u)=R(m_v\ot m_u)$ gives us
\begin{equation}\label{uns}
\sigma(c_{iv}\ot \overline{c_{ju}})=x_{uv}^{ji}
\end{equation}
for all $i$, $j$, $u$, $v=1,\cdots, n$. Hence, the equation (\ref{uns})
ensure the uniqueness of $\sigma$.

Now we shall prove the existence of $\sigma$. First we define
$$\sigma_0:{\cal M}^n(k)\ot {\cal M}^n(k)\to k,\qquad
\sigma_0(c_{iv}\ot {c_{ju}})=x_{uv}^{ji}
$$
for all $i$, $j$, $u$, $v=1,\cdots, n$. In order to
prove that $\sigma_0$ factorizes to a map
$\sigma :{\cal M}^n(k)\ot {\cal M}^n(k)/I(R)\to k$, we have
to show that $\sigma_0({\cal M}^n(k)\ot I(R))=0$.
For  $i$, $j$, $k$, $l$, $p$, $q=1,\cdots, n$, we have:
\begin{eqnarray*}
\sigma_0(c_{pq}\ot o(i,j,k,l))&=&
\sum_{v}x_{kv}^{ji}\sigma_0(c_{pq}\ot c_{vl})-
\sum_{\alpha}x_{kl}^{j\alpha}\sigma_0(c_{pq}\ot c_{i\alpha})\\
&=&\sum_{v}x_{kv}^{ji}x_{lq}^{vp}-
\sum_{\alpha}x_{kl}^{j\alpha}x_{\alpha q}^{ip}\\
\text{(from (\ref{acum}))}
&=&0
\end{eqnarray*}
Hence we have constructed $\sigma :{\cal M}^n(k)\ot {\cal M}^n(k)/I(R)\to k$
such that $R=R_{\sigma}$. It remains to prove that $\sigma$ is a
${\cal D}$-map. Let $c=c_{ij}$, $\overline{d}=\overline{c_{pq}}$. We have:
$$
\sum \sigma(c_{(1)}\ot \overline{d})\:\overline{c_{(2)}}=
\sum_v\sigma(c_{iv}\ot \overline{c_{pq}})\:\overline{c_{vj}}=
\sum_vx_{qv}^{pi}\overline{c_{vj}}
$$
and
$$
\sum \sigma(c_{(2)}\ot \overline{d})\:\overline{c_{(1)}}=
\sum_{\alpha}\sigma(c_{\alpha j}\ot \overline{c_{pq}})\:
\overline{c_{i\alpha}}=\sum_{\alpha} x_{qj}^{p\alpha}\overline{c_{i\alpha}}.
$$
Hence,
$$
\sum \sigma(c_{(1)}\ot \overline{d})\:\overline{c_{(2)}}-
\sum \sigma(c_{(2)}\ot \overline{d})\:\overline{c_{(1)}}=
\overline{o(i,p,q,j)}=0,
$$
i.e. $\sigma$ is a ${\cal D}$-map.

Suppose now that $R$ is bijective and let $S=R^{-1}$. Let $(y_{uv}^{ji})$
be a family of scalars of $k$ such that
$$
S(m_v\ot m_u)=\sum_{i,j}y_{uv}^{ji}m_i\ot m_j,
$$
for all $u$, $v=1,\cdots ,n$. As $RS=SR=Id_{M\ot M}$ we have
$$
\sum_{\alpha, \beta}x_{\beta \alpha}^{ip}y_{jq}^{\beta\alpha}=
\delta_{ij}\delta_{pq}, \qquad
\sum_{\alpha, \beta}y_{\beta \alpha}^{ip}x_{jq}^{\beta\alpha}=
\delta_{ij}\delta_{pq}
$$
for all $i$, $j$, $p$, $q=1,\cdots, n$. We define
$$
\sigma_{0}^{\prime}:{\cal M}^n(k)\ot {\cal M}^n(k)\to k, \quad
\sigma_{0}^{\prime}(c_{iv}\ot c_{ju}):=y_{uv}^{ji}
$$
for all $i$, $j$, $u$, $v=1,\cdots, n$. First we prove that
$\sigma_{0}^{\prime}$ is an inverse in the convolution algebra
$\Hom_k({\cal M}^n(k)\ot {\cal M}^n(k),\;k)$ of $\sigma_{0}$.
Let $p$, $q$, $i$, $j=1,\cdots, n$. We have:
\begin{eqnarray*}
\sum \sigma_0\Bigl((c_{pq})_{(1)}\ot (c_{ij})_{(1)}\Bigl)
\sigma_{0}^{\prime}\Bigl((c_{pq})_{(2)}\ot (c_{ij})_{(2)}\Bigl)&=&
\sum_{\alpha, \beta} \sigma_0(c_{p\alpha}\ot c_{i\beta})
\sigma_{0}^{\prime}(c_{\alpha q}\ot c_{\beta j})\\
&=&\sum_{\alpha, \beta}x_{\beta \alpha}^{ip}y_{jq}^{\beta\alpha}=
\delta_{ij}\delta_{pq}=\varepsilon(c_{ij})\varepsilon(c_{pq})
\end{eqnarray*}
and
\begin{eqnarray*}
\sum \sigma_0^{\prime}\Bigl((c_{pq})_{(1)}\ot (c_{ij})_{(1)}\Bigl)
\sigma_{0}\Bigl((c_{pq})_{(2)}\ot (c_{ij})_{(2)}\Bigl)&=&
\sum_{\alpha, \beta} \sigma_0^{\prime}(c_{p\alpha}\ot c_{i\beta})
\sigma_{0}(c_{\alpha q}\ot c_{\beta j})\\
&=&\sum_{\alpha, \beta}y_{\beta \alpha}^{ip}x_{jq}^{\beta\alpha}=
\delta_{ij}\delta_{pq}=\varepsilon(c_{ij})\varepsilon(c_{pq}).
\end{eqnarray*}
Hence, $\sigma_0\in \Hom_k({\cal M}^n(k)\ot {\cal M}^n(k),\; k)$ is
invertible in convolution. In order to prove that
$\sigma\in \Hom_k({\cal M}^n(k)\ot {\cal M}^n(k)/I(R),\; k)$ remains
invertible in the convolution, it is enought to prove that
$\sigma_{0}^{\prime}$ factorizes to a map
$\sigma^{\prime}:{\cal M}^n(k)\ot {\cal M}^n(k)/I(R)\to k$.
For $i$, $j$, $k$, $l$, $p$, $q=1,\cdots,n$ we have
$$
\sigma_{0}^{\prime}(c_{pq}\ot o(i,j,k,l))=
\sum_{v}x_{kv}^{ji}y_{lq}^{vp}-
\sum_{\alpha}x_{kl}^{j\alpha}y_{\alpha q}^{ip}
$$
As $S=R^{-1}$ and $R$ is a solution of the ${\cal D}$-equation we get that
$R^{23}S^{12}=S^{12}R^{23}$. Using the equation (\ref{unde}) we obtain
that $\sigma_{0}^{\prime}(c_{pq}\ot o(i,j,k,l))=0$. Hence,
$\sigma_{0}^{\prime}$ factorizes to a map
$\sigma^{\prime}:{\cal M}^n(k)\ot {\cal M}^n(k)/I(R)\to k$.

2. Suppose now that $R\tau=\tau R$. It follows that
\begin{equation}\label{22}
x_{uv}^{ji}=x_{vu}^{ij}
\end{equation}
for all $i$, $j$, $u$, $v=1,\cdots,n$. Let
$$C(R):={\cal M}^n(k)/I(R)$$
The rest of the proof is similar to the one in item 1.
The only thing we have to prove is that the map
$\sigma :{\cal M}^n(k)\ot {\cal M}^n(k)/I(R)\to k$ factorizes to a map
$\tilde{\sigma} :C(R)\ot C(R)\to k$.
We have
\begin{eqnarray*}
\sigma(o(i,j,k,l)\ot \overline{c_{pq}})&=&
\sum_v x_{kv}^{ji}\sigma(c_{vl}\ot \overline{c_{pq}})-
\sum_{\alpha} x_{kl}^{j\alpha}\sigma (c_{i\alpha}\ot \overline{c_{pq}})\\
&=&\sum_v x_{kv}^{ji}x_{ql}^{pv}-
\sum_{\alpha} x_{kl}^{j\alpha}x_{q\alpha}^{pi} \\
\text{(\mbox{using} (\ref{22}) )}
&=&\sum_v x_{kv}^{ji}x_{lq}^{vp}-
\sum_{\alpha} x_{kl}^{j\alpha}x_{\alpha q}^{ip}=0
\end{eqnarray*}
for all $i$, $j$, $k$, $l$, $p$, $q=1,\cdots,n$. Hence,
$\sigma (I(R)\ot {\cal M}^n(k)/I(R))=0$, i.e. $\sigma$ factorizes to a map
$\tilde{\sigma} :C(R)\ot C(R)\to k$.
\end{proof}

Applying the above theorem for the solution of the ${\cal D}$-equation
mentioned
in the previous section, we shall construct the coresponding ${\cal D}$-map
defined for the coalgebra $C={\cal M}^4(k)$.

\begin{examples}
1. Let $C={\cal M}^4(k)$ and $I$ be the two dimensional $k$-subsapce of $C$
with $\{c_{21}, c_{22}-c_{11} \}$ a $k$-basis. Then $I$ is a coideal of $C$.
Let $a$, $b$, $c$ scalars of $k$ and $(x_{ij}^{uv})$ given by the formulas
(\ref{scalari}). Then
$$
\sigma :{\cal M}^4(k)\ot {\cal M}^4(k)/I\to k,\quad
\sigma(c_{iv}\ot \overline{c_{ju}})=x_{uv}^{ji}
$$
for all $i$, $j$, $u$, $v=1,\cdots, 4$ is a ${\cal D}$-map.

2. Let $C={\cal M}^4(k)$, $q\in k$ and $I$ be the two dimensional coideal
of $C$ with $\{c_{21}, c_{12}+qc_{22}-qc_{11} \}$ a $k$-basis. Let
$(x_{ij}^{uv})$ be the scalars of $k$ of which
$$
x_{21}^{11}=-q,\quad x_{21}^{21}=1,\quad
x_{22}^{11}=-q^2,\quad x_{22}^{21}=q.
$$
and all others are zero.
Then
$$
\sigma :{\cal M}^4(k)\ot {\cal M}^4(k)/I\to k,\quad
\sigma(c_{iv}\ot \overline{c_{ju}})=x_{uv}^{ji}
$$
for all $i$, $j$, $u$, $v=1,\cdots, 4$ is a ${\cal D}$-map.

3. Let $C={\cal M}^4(k)$  and $I$ be the two dimensional coideal
of $C$ with $\{c_{12}, c_{21} \}$ a $k$-basis. Let $a$, $b\in k$ and
$$
\sigma :{\cal M}^4(k)\ot {\cal M}^4(k)/I\to k
$$
such that
$$
\sigma(c_{11}\ot \overline{c_{11}})=a,\quad
\sigma(c_{11}\ot \overline{c_{22}})=b
$$
and all others are zero. Then $\sigma$ is a ${\cal D}$-map.
\end{examples}

\end{document}